\documentclass{amsart}
\usepackage{lipsum}
\usepackage{amsmath}
\usepackage{amsthm}
\usepackage{amsfonts}
\usepackage{amssymb}
\usepackage{latexsym} 
\usepackage{graphicx}
\usepackage{mathtools}
\usepackage{color}
\usepackage{appendix}

\usepackage{amsmath,amsthm,amssymb,mathrsfs} 
\usepackage[T1]{fontenc}
\usepackage{ae,aecompl}
\usepackage{times}
\usepackage[english]{babel}
\usepackage[colorlinks,pagebackref,hypertexnames=false]{hyperref}
\usepackage[alphabetic,backrefs]{amsrefs}
\usepackage{bookmark}
\usepackage{microtype}

\usepackage[matrix,tips,graph,curve]{xy}

\usepackage{tikz}
\usetikzlibrary{decorations.markings}



\makeatletter 
\@addtoreset{equation}{section}
\makeatother

\theoremstyle{plain}
\newtheorem{theorem}[equation]{Theorem}

\newtheorem{lemma}[equation]{Lemma}

\theoremstyle{definition}

\newtheorem{remark}[equation]{Remark}

\newcommand{\IR}{\mathbb{R}}

\newcommand\reals{{\mathbb R}}

\renewcommand\Im{\mathrm Im}

\newcommand{\<}{\langle}
\renewcommand{\>}{\rangle}

\newcommand{\ra}{\rightarrow}                   

\def\d/{/\mspace{-6.0mu}/}

\def\reals{{\mathbb R}}
\def\cx{{\mathbb C}}

\def\Im{\,\mathrm{Im}\,}

\def\supp{\mathrm{supp}\,}

\def\O{{\mathcal O}}
\def\SS{{\mathbb S}}
\def\s{{\mathcal S}}

\def\phi{\varphi}

\def\be{\begin{eqnarray*}}
\def\ee{\end{eqnarray*}}
\def\ben{\begin{eqnarray}}
\def\een{\end{eqnarray}}

\def\lll{\left\langle}
\def\rrr{\right\rangle}

\def\L2R{L_{\text{Rest}}^2}

\def\11{\mathds{1}}

\def\tpsi{\tilde{\psi}}
\def\tchi{\tilde{\chi}}
\def\L2c{L^2_{\text{comp}}}

\def\th{\tilde{h}}

\def\uhi{u_{\text{hi}}}
\def\ulo{u_{\text{lo}}}

\newcommand{\p}{\partial}

\begin{document}

\setcounter{section}{0}

\title[Local smoothing]{Local smoothing for the Schr\"odinger equation
  on a multi-warped
  product manifold with inflection-transmission trapping}
\author{Hans Christianson and Derrick Nowak}

\date{\today}

\begin{abstract}
Geodesic trapping is an obstruction to dispersive estimates for
solutions to the
Schr\"odinger equation.    Surprisingly little is known
about solutions to the Schr\"odinger equation on manifolds with
degenerate trapping, since the conditions for degenerate trapping are
not stable under perturbations.  In this paper we extend some of the 
results of \cite{ChMe-lsm} on inflection-transmission type trapping on
warped product manifolds to the case of {\it multi}-warped products.
The main result is that the trapping on one cross section does not
interact with the trapping on other cross sections provided the
manifold has only one infinite end and only inflection-transmission
type trapping.

\end{abstract}

\maketitle

\setcounter{section}{0}

\section{Introduction}

In this paper, we study the effects of inflection-transmission type trapping on local smoothing 
estimates for solutions to the Schr\"odinger equation on a
multi-warped product manifold.
Inflection-transmission trapping on a warped product manifold was introduced in \cite{ChMe-lsm} by
Christianson-Metcalfe as a semi-stable type of trapping.
The warped product structure allows the authors to separate variables
and study an essentially one-dimensional problem.  
The purpose of this paper is to continue that study into the context of
a multi-warped product manifold where the trapping can occur on
different cross sections.  This breaks the symmetry of the single
warped product manifold so that the problem is no longer a
one-dimensional problem.

\subsection{Multi-warped product manifold}
The most familiar example of a {\it warped} product manifold is a
surface of revolution, which involves a defining curve revolved around
a line. This means the defining curve is warping the
circle at each point to change  the radius along the surface.    
The second most familiar warped product manifold is 
$\reals^n$ in polar coordinates.  That is, $\reals^n = \reals_+ \times \SS^{n-1}$
together with the metric
\[
g = dx^2 + x^2 g_{\SS^{n-1}}.
\]
Here we refer to $A(x) = x$ as the ``warping'' function.   Let $A(x) :
\reals_+ \to \reals$ be a smooth function satisfying $A(x)>0$
for $x>0$ and $A(x) 
\sim x$ near $x = 0$ and outside a compact set.  
Let 
$M$ be a compact
Riemannian manifold without boundary.  Then $X = \reals_+ 
\times M$
with the metric
\[
g = dx^2 + A^2(x) g_{M}
\]
is called a warped product with cross section $M$ and warping function
$A(x)$.  It is ``Euclidean'' outside a compact set because $A(x) = x$
outside a compact set, and it has one ``infinite end'' since we are
only working with $x \in \reals_+$ and $A(x) = x$ near $x = 0$.

A multi-warped product is a product of two or more cross section manifolds warped by
different warping functions.  We will assume our manifold is Euclidean
outside a compact set so that infinity looks like a compact product manifold
warped in the usual polar coordinates.  In this paper, we will
specialize to the case with only one infinite end.

A multi-warped product manifold is defined as follows:  Let $M_1, M_2,
\ldots M_N$ be compact Riemannian manifolds without boundary.  Denote the
corresponding metrics $g_{M_1}, \ldots , g_{M_N}$, and suppose they
have dimensions $n_1, \ldots , n_N$ respectively.  Let $A_1, \ldots,
A_N: \reals_+ \to \reals$ satisfy $A_j(x) >0$, $A_j(x) = x$ near $x=0$ and outside
a compact set.  Let 
\[
X = \reals_+ \times M_1 \times M_2 \times \cdots \times M_N
\]
with the metric
\[
g = dx^2 + A_1(x)^2 g_{M_1} + \ldots + A_N(x)^2 g_{M_N}
\]
Then X is a multi-warped product manifold with cross sections $M_1 , \ldots
M_N$.  It is Euclidean at infinity, since the metric is
\[
g= dx^2 + x^2(g_{M_1} + \ldots + g_{M_N})
\]
for $x$ outside a compact set.  The metric $g$ takes the same form in
a neighborhood of $x = 0$, so $X$ is Euclidean near $0$ as well.
Observe that the dimension of $X$ is $n_1 + n_2 + \ldots + n_N +1$.

Many of these assumptions about the geometry can be relaxed in various
ways without significantly changing the analysis in this paper.  It is also possible to study multi-warped product manifolds with two
ends, which just means the $A_j(x)$ are positive functions on $\reals$
which equal $|x|$ outside a compact set.  We will study the 
Schr\"odinger equation on such manifolds in a subsequent paper.

\section{Statement of Results}

Let $X$ be a Riemannian manifold with metric $g$, and let $-\Delta_g$
denote the corresponding Laplace-Beltrami operator.  
The Schr\"odinger equation on $X$ is
\begin{align}
  \begin{cases}
    \label{E:Sch}
  (D_t -\Delta_g) u(t,x) = 0 \text{ on } \reals_t \times X,
  \\
  u(0,x) = u_0(x),
\end{cases}
\end{align}
where $u_0$ is in some reasonable Sobolev space.  Here we use the
convention $D_t = \frac{1}{i} \p_t$.  Our goal is to understand how
the geometry of $X$ affects solutions to \eqref{E:Sch}.  In the
following subsection we construct a multi-warped product manifold with
inflection-transmission type trapping.

\subsection{Construction of the Manifold}

In order to make the present paper as clear as possible, we specialize
to the case where there are only two cross sections, both circles.

\begin{figure}[h]
  \caption{\label{F:Ajs} The functions $A_1$ and $A_2$.}
\centering
\begin{tikzpicture}
\draw[line width=0.5mm,->] (0,0) -- (12,0);
\draw[line width=0.5mm,->] (0,0) -- (0,10);
\draw[line width=0.3mm, dashed] (4,2) parabola (2,1);
\draw[line width=0.3mm, dashed] (4,2) parabola (6.9,3.67);
\draw[line width=0.3mm, dashed] (6.9,3.67) -- (10, 7);

\draw[fill] (4,2) circle [radius=0.05];
\node [below] at (4,2) {\textbf{$x=1$}};

\draw[fill] (8,6) circle [radius=0.05];
\node [above] at (8,6) {\textbf{$x=2$}};

\node [above] at (4,3.5) {\textbf{$A^2_2(x)$}};
\node [below] at (8,4.5) {\textbf{$A^2_1(x)$}};

\draw[line width=0.3mm, domain=2:5,smooth,variable=\y]  plot ({\y},{39.3/(4-(1/4)*\y)-10.2});
\draw[line width=0.3mm] (8, 6) parabola (5,4.1);
\draw[line width=0.3mm] (8, 6) parabola (10,7);
\draw[line width=0.4mm, black] (0,0) parabola (2.15,1.15);
\draw[line width=0.3mm, dashed] (4,2) parabola (2,1);
\begin{scope}
	\clip (11.67,8) circle (2cm);
	\draw[line width=0.4mm, black] (7.67, 5.94) parabola (11.67,9.04);
\end{scope}
\draw[fill] (8,6) circle [radius=0.05];
\node [above] at (8,6) {$x=2$};
\end{tikzpicture}
\end{figure}
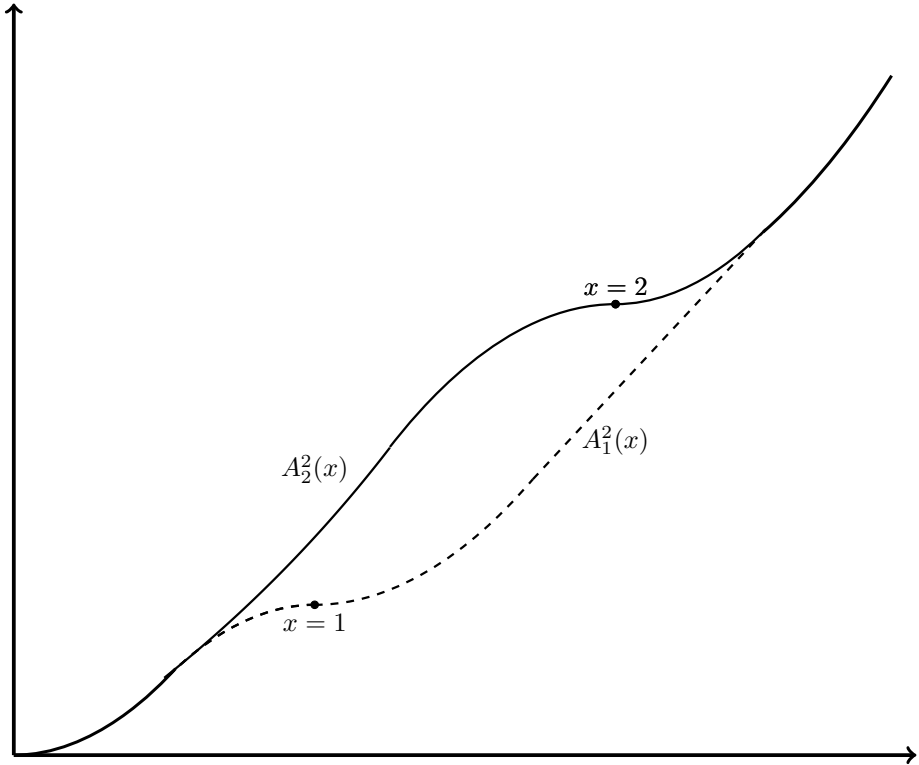

We  consider smooth functions $A_1,A_2$ and constants
$C_1,C_2,C_3,C_4,C_5,C_6,C_7,C_8$ such that for $j=1,2$, $A_j(x) = x$ for $x$ near
$0$ and outside a compact set, $A_j(x)>0$ for $x >0$, $A_j'(x) \geq 0$,
$$A_1^2(x) = \begin{cases}
C_1(x-1)^{2m_1+1}+C_2,& x\sim 1\\
\frac{1}{C_3-C_4x}, & x\sim 2
\end{cases}
$$
and
$$A_2^2(x) = \begin{cases}
\frac{1}{C5-C_6x}, & x\sim 1 \\
C_7(x-2)^{2m_2+1}+C_8,& x\sim 2
\end{cases}$$
where $A_1'(x)=0$ if and only if $x=1$ and  $A_2'(x)=0$ if and only if $x=2$.
Here $m_1$ and $m_2$ are positive integers.  The constants are needed to make sure such functions exist while
maintaining that $A_2^2,A_1^2$ have only one point where the
derivative is $0$.  We are also assuming that $A_1^{-2}(x)$ is linear
and decreasing
near $x=2$ and $A_2^{-2}(x)$ is linear and decreasing near $x=1$.
A sketch of $A_1$ and $A_2$ are found in Figure \ref{F:Ajs}.

Now let $X = \reals_+ \times \SS^1 \times \SS^1$ be a half line
crossed with two circles.  Let $\theta$ and $\omega$ parametrize the
circles, and let
\[
g = dx^2 + A_1^2(x) d \theta^2 + A_2^2(x) d \omega^2,
\]
making $X$ a multi-warped product manifold.

\begin{theorem}
  \label{T:main}
Let $(X,g)$ be the multi-warped product constructed above.  Suppose
$u$ solves \eqref{E:Sch} on $X$ with $u_0 \in \s(X)$.  Let $m = \max
(m_1, m_2)$.  Then  for each
$T>0$ there exists a constant $C$ such that
\begin{align}
  \label{E:main}
  \int_0^T \| \lll x \rrr^{-3/2} u \|^2_{H^1(X)} dt \leq C \| u_0
  \|^2_{H^{\frac{2m+1}{2m+3}}(X)}
  \end{align}

\end{theorem}

\begin{remark}
Due to the quasimode construction in \cite{ChMe-lsm} the estimate in \eqref{E:main} is sharp.
\end{remark}

 \begin{remark}
  The estimate \eqref{E:main} expresses that locally in space and on
  average in time the solution $u$ is $2/(2m+3)$ derivatives smoother
  than the initial data.  Because of this, estimate \eqref{E:main} is
  called a local smoothing estimate.  See Subsection \ref{SS:history} for
  motivation and history of local smoothing type estimates.
  \end{remark}

\begin{remark}
  Additionally, we state the result in the form listed because it follows naturally from the estimates below. However, since $u_0\in\s(X)$ we can commute with an appropriate psuedo-differential operator to get that for any $s,$
  $$\int_0^T||\<x\>^{-3/2}u||_{H^{s+\frac{2}{2m+3}}(X)}^2dt\leq C||u_0||_{H^{s}(X)}^2$$
\end{remark}

\begin{remark}
The power $-3/2$ in the weight function is not optimal, but helps our
computations later.
 We have assumed $u_0 \in\s$ to avoid any regularity issues, but a
  density argument can be used to extend this result to rougher
  initial data. 
\end{remark}

\begin{remark}
We again want  to emphasize that $A_1^2$ has an
inflection point of order $2m_1+1$ at $x=1$, $A_2^2$ has an inflection
point of order $2m_2+1$ at $x=2$ and that $A_1$ and $A_2$ give the
Euclidean metric near $x=0$ and when $x$ is large. We also make
$A_1^{-2}$ linear near $x=2$ and $A_2^{-2}$ linear near $x=1$ to make
some of the computation easier. However, we expect that this
conditioned  can be loosened and still give the same result.
\end{remark}

\subsection{ Motivation and History}
\label{SS:history}

The Schr\"odinger equation is one of a large family of {\it  dispersive} equations, which are equations whose solutions propagate
in a way that depends on the frequency of oscillation.  Dispersive
equations have conserved quantities, often expressing that the mass or size
of oscillations are preserved in time.  For the Schr\"odinger
equation on $\reals^n$,  the $H^s$ norm of a solution is preserved in
time.  In other words, at any time $t$, the solution has the same
regularity as the initial data.  
The local
smoothing effect for solutions to the Schr\"odinger equation expresses that, even
though a solution to the Schr\"odinger equation has  the same
regularity as the initial data, on average in time and locally in
space the solution is $1/2$ derivative smoother.

The local smoothing estimate for solutions to the Schr\"{o}dinger
equation on $\reals^n$ is that for any $T$ and any $\varepsilon>0$, there exists a $C>0$ such that 
$$\int_0^T \|\<x\>^{-1/2-\varepsilon}
e^{it\Delta}u_0\|_{H^{1/2}}^2dt\leq C \|u_0\|_{L^2}^2.$$
This type of estimate has been studied in a number of different
contexts with dispersive equations of varying orders
\cite{Sl,ConSau,Vega}.  These studies were extended to the case of
non-trapping asymptotically Euclidean manifolds in
\cite{CKS,MR1373768}.  That trapping necessarily causes a loss in
regularity was proved by Doi \cite{Doi2}.

There have been a number of results about manifolds with trapping.  If the trapping is unstable and non-degenerate, the loss in regularity
is logarithmic \cite{Bur-sm,Chr-NC,Chr-disp-1,Chr-QMNC,Dat-sm}.
Non-degenerate trapping allows the use of quantum Birkhoff normal
forms to have an invariant definition of hyperbolic trapping.  If the
trapping is unstable but degenerate, normal forms are not available so
the examples are limited.  In \cite{ChWu-lsm} the authors show there is a
local smoothing estimate with sharp polynomial loss.  In
\cite{ChMe-lsm} the authors introduce the semi-stable
inflection-transmission trapping, further studied in the present
paper, and demonstrate a local smoothing estimate with sharp
polynomial loss.  In \cite{Chr-inf-deg}, the author proves that
unstable but infinitely degenerate trapping causes a complete loss.

The intuition behind the non-trapping estimates is as follows: In $\reals^n$, if $u_0$ is sufficiently smooth, we can use the Fourier
transform to write down the solution:
\[
u(t,x) = c_n \int u_0(y) e^{i(-t|\xi|^2 + \xi\cdot (x-y))} dy d \xi,
\]
where $c_n$ is a dimensional constant.
Restricting our attention to  $\reals^2$, the solution
has phase function $-t \xi^2 + \xi (x-y)$ which is stationary when $-2t
\xi + (x-y) = 0$, or $x = y + 2t\xi$.  This means that a  solution at frequency $\xi$ propagates at
speed $2 \xi$.  This has the effect that a solution leaves a compact
set in space in time $t  \sim \xi^{-1}$.  Then integrating the $H^s ( \reals^2)$
norm in time gains $ \xi^{-1}$ over $|\xi|^{2s}|\hat{u} |^2$, or $1/2$
derivative on each copy of the solution $u$.

We also see from this heuristic that solutions propagate along geodesics in the sense
that they follow straight lines as they propagate out to infinity.  The
same is true on manifolds, as long as all geodesics go to infinity.
This is why trapping plays such an important role in local smoothing
estimates.  
When trapping occurs, wave packets can stay  coherent near the trapping
which means that our $\reals^2$ heuristic does not work any more, and
we expect some loss in regularity.


\subsection{Overview}

On a warped product manifold $X = \reals_+ \times M$ with metric $g =
dx^2 + A^2(x) g_M$,
the Laplacian is, up to lower order terms,
\[
- \Delta = -\p_x^2 - A^{-2}(x) \Delta_{g_M}.
\]
Let $\{ \phi_j(\omega) \}$ be the orthonormal basis of $L^2(M)$
consisting of eigenfunctions:
\[
-\Delta_{g_M} \phi_j = \lambda_j^2 \phi_j.
\]
 Then if $f:X \to \cx$ is sufficiently smooth, we can separate variables:
\[
f(x, \omega) = \sum f_j(x) \phi_j ( \omega),
\]
so that, up to lower order terms, 
\[
-\Delta f = \sum (-f_j'' + A^{-2}(x) \lambda_j^2f_j) \phi_j.
\]
On each eigenspace then one considers the operator $-\p_x^2 +
\lambda_j^2 A^{-2}(x).$  Rescaling $h = \lambda_j^{-1}$, we are led to
consider the operator $P = -h^2 \p_x^2 + V(x)$, where $V(x) =
A^{-2}(x)$.  The corresponding (semi-classical) symbol is $p = \xi^2 +
V(x)$.  In this reduced geometry, the replacement for the geodesic
flow is the Hamiltonian flow, and solutions propagate along this
flow.  The Hamiltonian system for this symbol is then
\[
\begin{cases}
  \dot{x} = 2 \xi, \\
  \dot{\xi} = -V'(x),
  \\
  x(0) = x_0, \\
  \xi(0) = \xi_0.
\end{cases}
\]
If  $V'(x_0) = 0$, then $(x,
\xi) = (x_0, \xi_0)$ is a ``trapped'' solution.  This corresponds to a
longitudinal periodic geodesic on the original warped product.

The question of local smoothing with loss then boils down to
understanding what happens to solutions of the one-dimensional
semi-classical problem  near critical points in phase space.  This
necessitates use of second microlocalization to get sharp estimates.
This analysis was done in the papers \cite{ChWu-lsm} with degenerate
unstable trapping, \cite{ChMe-lsm} for inflection-transmission type
trapping, and in \cite{Chr-inf-deg} for infinitely degenerate critical points.
The present paper is a continuation of this series of papers.

The motivation is to see how different kinds of trapping interact at
different frequencies in a relatively simple geometric setting.
 Our main result, however, is that
the trapped sets on each cross section do not see each other, so the loss in
local smoothing is the same as in \cite{ChMe-lsm}.

Nevertheless, there are a number of  things to prove.
Having a product of two compact manifolds as cross sections, one can
separate variables on both cross sections.  Then one is led to study a
one-dimensional problem with two frequency parameters.  This appears to be a
complicated mess comparing different frequencies.  However, we can
separate variables in one cross section alone, which leaves us with a
two-dimensional problem with one parameter.
Since we are only
separating variables in one direction, we do have to deal with
derivatives in the other direction.  However, a detailed microlocal
frequency localization allows us to handle this problem.  The fact
that the trapping on one cross section does not see the trapping on
the other cross section is special to the one ended case and not
expected to hold in general.




\section{Local Smoothing away from the trapping}

Now that we have $A_1$ and $A_2$ defined, consider the product manifold $\IR_+\times \SS^1\times \SS^1$ with the metric
$$g=dx^2+A_1(x)^2d\theta^2+A_2(x)^2 d\omega^2.$$
Then, the laplacian is given by
$$\Delta_g = \p_x^2+A_1(x)^{-2}\p_{\theta}^2+A_2(x)^{-2}\p_{\omega}^2+(A'_1(x)A_1^{-1}(x)+A'_2(x)A_2^{-1}(x))\p_x$$
Next we use a transformation to get rid of the $\p_x$
term. Consider the unitary tranformation $T:L^2 (X, dV_g) \to L^2(X,
dx d \theta d \omega)$ given by 
$$Tu = A_1^{1/2}(x)A_2^{1/2}(x)u$$ and set
$$\tilde{\Delta}=T^{}\Delta_g T^{-1}.$$
This gives 
$$\tilde{\Delta} = \p_x^2+A_1^{-2}\p_{\theta}^2+A_2^{-2}(x)\p_{\omega}^2+ V(x)$$
where
\begin{align*}
V=&\frac{1}{4}A_1'(x)^2A_1(x)^{-2}-\frac{1}{2}A_1''(x)A_1^{-1}(x)\\
&+\frac{1}{4}A_2'(x)^2A_2(x)^{-2}-\frac{1}{2}A_2''(x)A_2^{-1}(x)\\
&-\frac{1}{2} A_1(x)^{-1}A_2(x)^{-1}A_1'(x)A_2'(x).
\end{align*}
This $V$ is similar to the single warped product case except we have a cross term of $$\frac{1}{2} A_1(x)^{-1}A_2(x)^{-1}A_1'(x)A_2'(x).$$

Next we want to do a positive commutator argument to get local smoothing away from $x=1,x=2$.
Let $u$ be a solution to $(D_t-\tilde{\Delta})u=0.$
Notice that $\tilde{\Delta}$ is of a similar form to \cite{ChWu-lsm}.
Let us take $B= f(x)\p_x$ for some general $f\in C^2(\IR)$ such that $f,f', f''$ are all bounded and then we will reduce to a specifc case.
\begin{align*}
[\tilde{\Delta},B]=& 2f'(x) \p_x^2+f''(x)\p_x\\
	&+2 A_1'A_1^{-3}f(x) \p_{\theta}^2+ 2A_2'A_2^{-3} f(x) \p_{\omega}^2+V'(x)f(x)
\end{align*}

\begin{remark}
Note that
$$\<u,v\>= \int_{\IR_+}\int_{\SS^1}\int_{\SS^1} u\bar{v} dx d\theta d\omega.$$
\end{remark}
and that
$$iB-(iB)^*=i[f(x),\p_x].$$ Hence,

\begin{align*}
0 &= \int_0^T\int_{\IR_+\times \SS^1\times \SS^1} u \overline{(f(x)D_x(D_t-\tilde{\Delta})u)}dx d\theta d\omega dt\\
	&= \int_0^T\int f(x)(D_x u)\overline{((D_t-\tilde{\Delta})u)}dx d\theta d\omega dt \\
	&+ \int_0^T \int (iB-(iB)^*)u\overline{((D_t-\tilde{\Delta})u)}dxd\theta d\omega dt\\
	&= i\<f(x)D_xu,u\>|_0^T+\int_0^T\<(D_t-\tilde{\Delta})i^{-1}Bu,u\>dt.
\end{align*}
This follows from integrating $\p_t(\<f(x)D_x u,u\>)$ in $t$, using $D_t u = \tilde{\Delta} u$ and integrating by parts. It is the same computation as \cite{ChWu-lsm} and the next step in the paper follows through as well.
Using the notation that $P=D_t-\tilde{\Delta}$,
\begin{align*}
0 &= 2i \,\Im\int_0^T\<i^{-1}BPu,u\> dt\\
	&=\int_0^T \<i^{-1}BPu,u\> dt - \int_0^T\<u, i^{-1}BPu\>dt \\
	&=\int_0^T\<[i^{-1}B,P],u,u\>dt - i \<f(x)D_x u,u\>|_0^T.
\end{align*}
Since $B$ is independent of $t$ this gives
$$\int_0^T \<[B,-\tilde{\Delta}]u,u\> dt = -\<f(x)D_x u, u \>|_0^T.$$

Let us reduce to the specific case of a function $f(x)$ given by
the following: let $\zeta(x)$ be a smooth function satisfying
$\zeta(x) \equiv 1$ near $x = 0$, $\zeta(x) >0$ for all $x$, and $|
\zeta(x) | \sim \lll x \rrr^{-3}$ for large $x$.  Such a
$\zeta$ is integrable, so let
\[
f(x) = \int_0^x \zeta(t) dt.
\]
Then $f(x) = x$ near $x = 0$, and there exists a constant $c>0$ such
that $f'(x) \geq c \lll x \rrr^{-3}$ for $x \geq 0$.  The power $-3$
here is much bigger than needed, but we have chosen it so that our
computation are easier.  We simply are matching the power of 
each $A_j^{-3} \sim x^{-3}$ as $x \to \infty$.  

  The restriction that $f(x)$ is linear near $x = 0$ is just to
maintain all the properties of Euclidean polar coordinates near $x=0$.   
Integrating by parts yields
\begin{align*}
&\int_0^T -\<2f'(x)\p_xu, \p_x u\>-\<2A_1'A_1^{-3}f(x) \p_\theta u,
  \p_\theta u\>-\<2A_2'A_2^{-3}f(x) \p_\omega u, \p_\omega u\> dt\\
&=-\<f(x) D_x u, u\>|_0^T+\int_0^T\<f''(x)\p_x u,u \>-\<V'(x)f(x)u,u\>dt
\end{align*}
Let us quickly remark again that, since each $A_j(x) = x$ for $x$ near
$0$ and $f(x) = x$ for $x$ near $0$, we have $A_j^{-3}(x) f(x) =
x^{-2}$ near $x = 0$.  We also have $V'(x) = 0$ near $x = 0$, so all
terms agree with the corresponding Euclidean terms near $x = 0$.  

Taking the absolute value of both sides and noting that $f'$, $A_1'A_1^{-3}f(x)$, and $A_2'A_2^{-3}f(x)\geq 0$ yields
\begin{align*}
&\int_0^T \|\sqrt{2f'(x)}\p_xu\|^2_{L^2}+
  \|\sqrt{2A_1'A_1^{-3}f(x)}\p_\theta u\|^2_{L^2}
  +\|\sqrt{2A_2'A_2^{-3}f(x)} \p_\omega u\|^2_{L^2}dt \\
&\leq C_1|\<D_x u, u\>|_0^T|+\int_0^TC_2|\<\p_x u,u \>|+C_3|\<u,u\>|dt
\end{align*}

Note that each term on the RHS is bounded by $C_T\|u_0\|_{H^{1/2}}$ for some constant $C_T.$
Next, we want to provide lower bounds on the $\sqrt{2f'(x)}$, $\sqrt{2A_1'A_1^{-3} f(x)}$, and $\sqrt{2 A_2'A_2^{-3}f(x)}$ terms.

First we want to bound the $\p_x^2$.  Note that $f'(x)= \zeta(x)$
defined above, so there exists a positive constant $c>0$ such that
$$\|\sqrt{2f'(x)}\p_xu\|^2_{L^2}\geq c\|\<x\>^{-3/2}\p_x u\|^2_{L^2}.$$

To get the correct lower bounds for the $\sqrt{2A_1'A_1^{-3} f(x)}$, and $\sqrt{2 A_2'A_2^{-3}f(x)}$ terms we will have to estimate $A_1'A_1^{-3}$ and $A_2'A_2^{-3}.$

\subsection{$A_1$ and $A_2$ estimates}

We have that
$$A_1^2(x) \sim \begin{cases}
C_1(x-1)^{2m_1+1}+C_2,& x\sim 1\\
x^2, &x \text{ away from 1}
\end{cases}
$$
So, near $x=0$
$$f(x)A_1'(x)A_1^{-3}(x)=\frac{1}{x^2}\geq C\frac{1}{x^2\<x\>^{1}}\geq C\frac{(x-1)^{2m_1}}{x^2\<x\>^{1+2m_1}}$$
Near $x=1$
$$f(x)A_1'(x)A_1^{-3}(x)\sim \frac{(x-1)^{2m_1}}{(1+(x-1)^{2m_1+1})^{3/2}}\geq C\frac{(x-1)^{2m_1}}{x^2\<x\>^{1+2m_1}}$$
When $x$ is large
$$f(x)A_1'(x)A_1^{-3}(x)\sim \frac{1}{x^3}\geq C\frac{1}{x^2\<x\>}
\geq C\frac{(x-1)^{2m_1}}{x^2\<x\>^{1+2m_1}}$$

Now just to be careful, we can consider compact sets $[\varepsilon,
  1-\varepsilon]$ and $[1+\varepsilon, K]$ for $K$ sufficiently large
and $\varepsilon$ small to handle the situation where we do not know the exact form of $A_1^2$. We know that on this region $f(x), A_1'(x)> 0$ so we can find $C>0$ sufficiently small so that
$$f(x)A_1'(x)A_1^{-3}(x)\geq C\frac{(x-1)^{2m_1}}{x^2\<x\>^{1+2m_1}}\, \, x\in [\varepsilon, 1-\varepsilon]\cup[1+\varepsilon, K]$$

With $A_2(x),$ the only difference is the inflection point is at
$x=2$ and we replace $m_1$ with $m_2$. This does not change the qualitative behavior of the
estimates. We just need estimates near $x=2$ instead of $x=1$ and we
will get $(x-2)^{2m_2}$ in the numerator instead of $(x-1)^{2m_2}$.  This
proves the following Lemma:
\begin{lemma}

  \label{L:initial}
Let $u$ be a solution to \eqref{E:Sch} on our manifold $X$ with
initial data $u_0 \in \s(X)$.  Then for each $T>0$, there exists a
constant $C>0$ such that
\begin{align}
  \int_0^T& ( \| \lll x \rrr^{-3/2} \p_x u\|^2   +  \|
  (x-1)^{m_1} \lll x \rrr^{-1/2-m_1}A_1^{-1} \p_\theta u \|^2\notag  \\
  & \quad + \|
(x-2)^{m_2} \lll x \rrr^{-1/2 -m_2} A_2^{-1} \p_\omega u \|^2 )dt
  \notag \\
& \leq
C \| u_0 \|^2_{H^{1/2}(X)}. \label{E:lsm-trapping}
\end{align}

\end{lemma}

\begin{remark}
The estimate \eqref{E:lsm-trapping} expresses that there is perfect
local smoothing in the radial $x$ direction with a loss at the trapped
set on each copy of $\SS^1$.  It is also clear that the statement of
Theorem \ref{T:main} could be sharpened to have loss only in $\theta$
and $\omega$ derivatives.  However, we have stated the theorem in the
simplest possible way to be clear.

  \end{remark}

\section{Separation of variables}
Consider the operator
$P_1=P_0+V(x)=-\p_x^2-V_1\p_\theta^2-V_2\p_\omega^2-V(x)$ where
$V_j=A_j^{-2}$ and $V(x)$ contains derivatives  of $A_j$ as shown above. Define a function $\varphi(x)\in C_c^\infty$ such that $0\leq \varphi(x)\leq1$, $\varphi(x)\equiv 1$, on $x\in [1-\varepsilon, 1+\varepsilon]$ for $1/4>\varepsilon>0$ and $\supp(\varphi)\subset [1-2\varepsilon, 1+2\varepsilon].$  Since we have local smoothing away from $x=1$ we can localize near this point. We do this now so that we can define a Fourier transform properly and do not have to worry about any integrability issues near $x=0$ due to the metric.

Now separate one variable at a time, starting with $\theta$.  Write
\[
u = \sum u_k (t, x , \omega) e^{i k \theta}, \,\, u_0 = \sum u_{0,k}
(x , \omega) e^{i k \theta}.
\]

Then each $u_k$ satisfies:
\[
(D_t + P_k-V) ((\phi u)_k) = 2\varphi'(x) \p_x u_k+\varphi''(x) u_k
\]
where
\[
P_k = -\p_x^2 + k^2 V_1 - V_2 \p_\omega^2.
\]
Note that $\varphi', \varphi''$ are compactly supported away from $x=0$ and $x=1$.

Below we will drop the subscript $k$ for notational purposes.
Now we want to decompose the frequency into high and low angular
frequency parts.  
The high frequency  part is when the frequency in the $\theta$ direction is large compared to the frequency in the $x$ direction.
Consider an even bump function $\psi\in C^\infty_C(\IR)$ which is $1$ for
$|r|\leq \varepsilon$ and vanishes for $|r|\leq 2\varepsilon$ for
$\varepsilon>0$ small. Define $$\uhi
= \psi(D_x/k) (\phi u),\, \ulo = (1 - \psi) (\phi u).$$ 
Since $\varphi$ provides a cutoff near $x=1$ and away from zero, we can define $\psi(D_x/k)$ in the usual way.

Now using the definition of $\ulo$ and the fact that $D_t u = -(P_k-V)u$ we get that.
\begin{align*}
  (D_t + P_k-V) \ulo & = [P_k -V ,(1-\psi)\varphi] u \\
  & = (1-\psi) [-\p_x^2 ,\varphi] u + [k^2V_1 - V_2
    \p_\omega^2 -V, -\psi ](\varphi  u)\\
  & = (1-\psi) (-2 \varphi' \p_x - \varphi'') u + [k^2 V_1-V , -
    \psi](\phi u) + [-V_2 \p_\omega^2 , - \psi ](\phi u)  \\
  & = (1-\psi) (-2 \varphi' \p_x - \varphi'') u + kL_1 (\phi u)
  -\frac{1}{k}L_2\p_\omega^2 (\phi u)
\end{align*}
Here $L_1$ and $L_2$ are semi-classical pseudo-differential operators
(with parameter $|k|^{-1}$)
of order zero with wavefront set contained in $\{ \psi' (\xi/k) \neq 0 \}
\subset \{ \varepsilon \leq |\xi|/|k| \leq 2 \varepsilon \}$, so we
observe 
$|D_x| \sim |k|$ on the wavefront set of $L_1$ and $L_2$.  We will use
this shortly.
Now combining the above statements gives
\[
(D_t + P_k-V ) \ulo = k L_1 (\phi u) - \frac{1}{k} L_2 \p_\omega^2 (\phi u)-(1-\psi)(2 \varphi' \p_xu+ \varphi'' u).
\]

We now run the commutator argument, but insert a cutoff $\chi_1(x)$
with $\chi_1 \equiv 1$ on $\supp ( \varphi)$ near $x = 1$ and $\chi_1 \equiv 0$ near $x =
2$.  Let us also assume that $\chi_1^{1/2}$ is still smooth.  Then with $B = f(x) \p_x$ as before, recalling that $f'(x) = \zeta(x)$, 
\begin{align}
  \int_0^T &  \lll \chi_1 [D_t + P_k-V , B ] \ulo, \ulo \rrr dt \label{E:comm-11}\\
  & = \int_0^T \lll \chi_1 ( -2\zeta(x) \p_x^2 - k^2 V_1' f (x)
  + V_2' f (x) \p_\omega^2) \ulo, \ulo \rrr dt \notag \\
  & + \int_0^T \lll \chi_1(- f''(x) \p_x+f V') \ulo, \ulo \rrr dt. \notag
\end{align}
The last line in \eqref{E:comm-11} has only one $x$ derivative, so is 
 bounded as follows:
\[
\left| \int_0^T \lll \chi_1(- f''(x) \p_x+f V') \ulo, \ulo \rrr dt
\right| \leq C_T \| u_0 \|^2_{H^{1/2}(X)}.
\]
The $\p_\omega^2$ term in the second line of \eqref{E:comm-11} is further estimated as follows: we know that
for $j = 1,2$, $V_j'(x)  \leq 0$ and our function $f \geq 0$, so 
\begin{align*}
  \int_0^T&  \lll \chi_1 V_2'f(x)\p_\omega^2 \ulo, \ulo \rrr dt \\
  & = -\int_0^T \lll \chi_1 V_2'f(x) \p_\omega \ulo, \p_\omega \ulo \rrr
  dt \\
  & \geq 0.
\end{align*}
We also know that $-\chi_1 f V_1' \geq 0$, so that
\begin{align} 
  \int_0^T & \lll \chi_1 ( -2\zeta(x) \p_x^2 - k^2 V_1' f (x)
  + V_2' f (x) \p_\omega^2) \ulo, \ulo \rrr dt \notag \\
  & \geq \int_0^T \lll -2 \chi_1  \zeta(x) \p_x^2 
   \ulo, \ulo \rrr dt.\label{E:lower-V2}
\end{align}

The next issue is to observe that $V_1'(1) = 0$, so does not help us
eliminate the vanishing at $x = 1$ in \eqref{E:lsm-trapping}.
However, we  observe that on the wavefront set of $\ulo$, we have $|k |
\lesssim |D_x|$, so we want to use the G{\aa}rding inequality to estimate
$k$ in terms of $D_x$.  Recall
that $\chi_1 V_1'f \leq 0$ and has compact support so 
 the G{\aa}rding inequality implies there exists a constant $C>0$ such that
\[
\lll k^2 \chi_1  \ulo, \ulo \rrr \leq -C \lll \zeta(x) \p_x^2
\ulo , \ulo \rrr + \O(1) \| \ulo \|^2_{H^{1/2}(X)}.
\]
Combining this with \eqref{E:comm-11} and \eqref{E:lower-V2}
\begin{align*}
  \int_0^T & \lll \chi_1 ( k^2 \ulo) , \ulo \rrr dt \\
  & \leq -C \lll \zeta(x) \p_x^2
  \ulo , \ulo \rrr + \O(1) \| \ulo \|^2_{H^{1/2}(X)} \\
  & \leq C  \int_0^T \lll \chi_1 ( -2\zeta(x) \p_x^2 - k^2 V_1' f (x)
  + V_2' f (x) \p_\omega^2) \ulo, \ulo \rrr dt + \O(1) \| \ulo
  \|^2_{H^{1/2} (X)} \\
  & = C 
 \int_0^T   \lll \chi_1 [D_t + P_k-V , B ] \ulo, \ulo \rrr dt + \O_T(1)
 \|u_0 \|^2_{H^{1/2}(X)}.
\end{align*}
Rearranging and using energy estimates, we have
\begin{equation}
\int_0^T \lll \chi_1 k \ulo, k \ulo \rrr dt \leq C_T \| u_0
\|^2_{H^{1/2}(X)} +C \left|\int_0^T   \lll \chi_1 [D_t + P_k-V , B ]
\ulo, \ulo \rrr dt\right|. \label{E:unpack-0}
\end{equation}

Now we unpack the commutator term $\int_0^T\< \chi_1[D_t+P_k-V,
  B]u_{lo},u_{lo}\>dt$. Integrating by parts yields,
\begin{align}
\left |\int_0^T\< \chi_1[D_t+P_k-V, B]u_{lo},u_{lo}\>dt\right| &\leq
2\left|\int_0^T \<\chi_1 B u_{lo},(D_t+P_k-V)u_{lo}\>dt\right| \label{E:unpack}\\
&+\left|\int_0^T \< B u_{lo}, 2\chi'_1\p_x u_{lo}\>dt\right| +
\left|\int_0^T\<Bu_{lo}, \chi''_1 u_{lo}\>dt\right|\notag \\
&+\left|\int_0^T\<(D_t+P_k-V)u_{lo}, (\chi_1 f)'u_{lo}\> dt\right|.
\notag 
\end{align}

We will examine each line of this estimate separately.  The key thing
to observe is that, since $B = f(x) \p_x$,  the first line in \eqref{E:unpack} has the highest
number of derivatives so will require the most work.  The terms with
just $\p_x$ derivatives can be controlled by our initial estimate in
Lemma \ref{L:initial}.  Hence,  
due to perfect local smoothing in the $x$ direction and energy
estimates, we can bound the two terms on the middle line of
\eqref{E:unpack}:
\begin{equation}
 \left|\int_0^T \< B u_{lo}, 2\chi'_1\p_x u_{lo}\>dt\right| +
\left|\int_0^T\<Bu_{lo}, \chi''_1 u_{lo}\>dt\right| \leq
C_T\|u_0\|^2_{H^{1/2}}. \label{E:unpack-middle}
\end{equation}

Now for the first and last line in \eqref{E:unpack} we want to use the fact that
\[
(D_t + P_k-V ) \ulo = k L_1 (\phi u) - \frac{1}{k} L_2 \p_\omega^2
(\phi u)-(1-\psi)(2 \varphi' \p_xu+ \varphi'' u).
\]

We can use the fact that $\varphi'$ and $\varphi''$ are compactly
supported away from $0$ and perfect local smoothing in the $x$
direction to get that 

\begin{align}\left|\int_0^T\<2(1-\psi) \varphi' \p_xu+(1-\psi) \varphi'' u,
(\chi_1 f )'u_{lo}\> dt\right|\leq
  C\|u_0\|_{H^{1/2}}^2 \label{E:unpack-2-1}
\end{align}
from the last line of \eqref{E:unpack},
and
\begin{align}
  \left|\int_0^T\<2(1-\psi) \varphi' \p_xu+(1-\psi) \varphi'' u, \chi
  B u_{lo}\> dt\right|\leq C\|u_0\|_{H^{1/2}}^2 \label{E:unpack-1-1}
  \end{align}
from the first line of \eqref{E:unpack}, 
for some constant $C$.

Next we want to handle the $k L_1 (\phi u) - \frac{1}{k} L_2 \p_\omega^2 (\phi u)$
term coming from the last line in \eqref{E:unpack}. To do this we can
use the fact that $\chi_1$ and $\chi'_1$ are supported away from $x=2$
so that we have perfect local smoothing in the $\omega$
direction according to Lemma \ref{L:initial}.  Hence
\begin{align}
  \label{E:unpack-2-2}
\left| \int_0^T \lll \frac{1}{k} L_2  \p_\omega^2 (\phi u) , (\chi_1 f)'\ulo \rrr dt
\right| \leq C_T\frac{1}{|k|} \| u_0 \|^2_{H^{1/2}}.
\end{align}
Now let $\tpsi$ be a smooth, even, compactly supported  bump function with $\tpsi (s) \equiv 1$ on
$\supp ( \psi (s) )$.  Let $\tchi_1$ be a smooth compactly supported
function such that $\tchi_1(s) \equiv 1$ on the support of $\chi_1$
but still supported away from $x = 0$ and $x = 2$.  Then
\[
\tpsi(D_x/k) L_1 \tpsi(D_x/k)  = L_1 + \O(|k|^{-\infty}), 
\]
which gives
\begin{align}
  & \left| \int_0^T \lll k L_1 (\phi u) , (\chi_1 f)' \ulo \rrr dt \right| \notag \\
  & \quad =  \left| \int_0^T \lll   \tchi_1 k\tpsi  L_1 \tpsi (\phi u) , (\chi_1 f)'  \ulo \rrr dt
  \right|  + C_T \| u_0 \|^2_{H^{1/2}}\notag \\
  & \quad =  \left| \int_0^T \lll  \tchi_1 k L_1 \tpsi (\phi u) , (\chi_1 f)'  \tpsi
  (\phi u) \rrr dt
  \right|  + C_T \| u_0 \|^2_{H^{1/2}}\notag \\
  & \quad \leq C \int_0^T |k| \| \tchi_1 \tpsi (\phi u) \|^2 dt + C_T \| u_0
  \|^2_{H^{1/2}}. \label{E:tpsi-1}
\end{align}

Combining \eqref{E:unpack-2-1} with \eqref{E:unpack-2-2} and
\eqref{E:tpsi-1}, we estimate the the last line in \eqref{E:unpack}:
\begin{align}
  & \left|\int_0^T\<(D_t+P_k-V)u_{lo}, (\chi_1 f)'u_{lo}\> dt\right|
  \notag \\
  & \quad \leq C \int_0^T |k| \| \tchi_1 \tpsi (\phi u) \|^2 dt + C_T \| u_0
  \|^2_{H^{1/2}}.  \label{E:unpack-2}
  \end{align}



We now proceed with the first line in \eqref{E:unpack}.  We have
already estimated the lowest order parts in \eqref{E:unpack-1-1}.  We
will deal with the term with $L_1$ last.  That means we need to
estimate
\[
\left| \int_0^T \lll \chi_1 B \ulo, \frac{1}{k} L_2 \p_\omega^2 (\phi u) \rrr
dt \right|.
\]
The difficulty is that there is one $x$ derivative in $B$ and two
$\omega$ derivatives.  We expect the $1/k$ to essentially remove one
derivative to use Lemma \ref{L:initial} away from $x = 1$.  However,
this requires some careful observations.

Since the
wavefront set of $L_2$ is contained where $\psi' \neq 0$, we have $|
D_x/k | \sim \varepsilon>0$ on the wavefront set of $\psi'$.
Recall that $\tchi_1$ is a bump function satisfying $\tchi_1 \equiv 1$ on
$\supp \chi_1$ but $\tchi_1 \equiv 0$ near $x = 2$.  We also choose a bump
function 
$\tpsi_1(r)$ satisfying $\tpsi_1(r) \equiv 1$ on $\supp \psi'(r)$ but $\tpsi_1(r)
\equiv 0$ near $r = 0$.  The point is that then $\left(\frac{\p_x}{ k}
\right)\tpsi_1(
D_x/k)$ is a bounded operator on $L^2$.

Then
\begin{align}
  \int & \lll \chi_1 B \ulo , k^{-1} L_2 \p_\omega^2 (\phi u) \rrr dt \notag \\
  & = \int \lll \chi_1 k^{-1} f(x) \p_x  \tpsi_1\left(\frac{D_x}{k}\right)  \p_\omega \ulo , \tchi_1
  L_2 \p_\omega u \rrr dt + O(k^{-\infty}) \int_0^T \|\tilde{\chi_1} \p_\omega
  u\|^2 dt \notag \\
  & \leq C\int \left( \left\| \chi_1(x) \left( \frac{\p_x}{k} \right) \tpsi_1\left(\frac{D_x}{k}\right) \p_\omega \ulo \right\|^2 + \| \tchi_1 L_2
  \p_\omega (\phi u) \|^2 \right) dt + O(k^{-\infty}) \int_0^T \|\tilde{\chi_1} \p_\omega
  u\|^2 dt . \label{E:unpack-1-2}
\end{align}
The operator
\[
\chi_1 (x) \left( \frac{ \p_x }{k} \right) \tpsi(D_x/k)
\]
is bounded on $L^2$ and supported away from $x = 2$.  Similarly, the
operator $\tchi_1 L_2$ is bounded on $L^2$ and supported away from $x
= 2$.  The $\O( |k|^{-\infty})$ term has $\tchi_1  \p_\omega u$, which
is again supported away from $x = 2$.
Lemma
\ref{L:initial} guarantees perfect local smoothing in the $\omega$
direction away from $x = 2$, so applying Lemma \ref{L:initial} to
\eqref{E:unpack-1-2} yields
\begin{align}
\left| \int  \lll \chi_1 B \ulo , k^{-1} L_2 \p_\omega^2 (\phi u) \rrr dt
\right| \leq C \| u_0 \|^2_{H^{1/2}}. \label{E:unpack-1}
  \end{align}


Combining \eqref{E:unpack-0} with \eqref{E:unpack-middle},
\eqref{E:unpack-2}, and \eqref{E:unpack-1}, we have
\begin{align}
  \int_0^T & \lll \chi_1 k^2 \ulo,  \ulo \rrr dt \leq C_T \| u_0
\|^2_{H^{1/2}(X)} +C \left|\int_0^T   \lll \chi_1 [D_t + P_k-V , B ]
\ulo, \ulo \rrr dt\right| \notag \\
& \leq C_T \| u_0 \|^2_{H^{1/2}} + C \left| \int_0^T \lll \chi_1 B
\ulo , k L_1 (\phi u) \rrr dt \right|\notag  \\
& \leq C_T \| u_0 \|^2_{H^{1/2}}  + C \left| \int_0^T \lll \chi_1 B
\ulo , \tchi_1 \tpsi k L_1 (\phi u) \rrr dt \right| \notag \\
& \leq C_T \| u_0 \|^2_{H^{1/2}} + C \int_0^T \| \tchi_1 \tpsi kL_1 (\phi u)
\|^2 dt \notag \\
& \leq C_T \| u_0 \|^2_{H^{1/2}} + C \int_0^T \| \tchi_1 \tpsi k(\phi u)
\|^2 dt.\label{E:low-freq}
\end{align}

Finally, we observe that, since $\uhi = \psi (\phi u)$, we have
\[
\| k \chi_1 \uhi \| \leq \| k \tchi_1 \uhi \| = \| k \tchi_1 \tpsi (\phi u)
\| + \O(1) \| u \| .
\]
That means
\begin{equation}
\int_0^T \| k \chi_1 \uhi \|^2 dt \leq C\int_0^T \| k \tchi_1 \tpsi (\phi u)
\|^2 dt + C_T \| u_0 \|^2_{H^{1/2}}. \label{E:hi-1}
\end{equation}
According to \eqref{E:low-freq}, we can estimate the low frequency
part of $u$ in terms of a quantity similar to the high frequency
estimate \eqref{E:hi-1}.
So for both $\uhi$ and $\ulo$, it suffices to estimate the high
frequency part:
\begin{align*}
  \int_0^T \| k \tchi_1 \tpsi (\phi u)
  \|^2 dt
\end{align*}
where $\tchi_1$ is supported near $x = 1$ and $\tpsi$ has compact
support.


\section{The high frequency estimate}

We use the $F F^*$ type argument employed in \cite{ChWu-lsm} and
\cite{ChMe-lsm}.  Let us drop the tilde notation and consider
functions $\chi(x)$ supported near $x = 1$ and supported away from $x
= 0$ and $x = 2$, as well as $\psi(D_x/k)$ micro-supported near $0$.  
Let $F(t)$ be defined by
$$F(t)g = \chi(x)\psi(D_x/k)k^r e^{-it(P_{k}-V)}g,$$
where $e^{-it(P_{k}-V)}$ is the free propagator. 
We want to show that for $r=\frac{2}{2m+3}$ we have a mapping $F:L_x^2\ra L^2([0,T])L_x^2,$ since then
$$\|k^{1-r}F(t)u_0\|_{L^2([0,T]);L^2)}\leq C\|k^{1-r}u_0\|_{L^2}$$
is the desired local smoothing estimate. We have such a mapping if and only if $FF^*:L^2([0,T])L_x^2\ra L^2([0,T])L_x^2$. We compute
$$FF^*f(x,t)=\psi(D_x/k)\chi(x)k^{2r}\int_0^T e^{i(t-s)(P_{k}-V)}\chi(x)\psi(D_x/k) f(x,s)ds,$$
and need to show that $\|FF^*f\|_{L^2L^2}\leq C\|f\|_{L^2L^2}$. Now write $FF^*f(x,t)=\psi\chi(v_1+v_2),$ where
$$v_1= k^{2r}\int_0^t e^{i(t-s)(P_{k}-V)}\chi(x)\psi(D_x/k)f(x,s)ds,$$
and
$$v_2= k^{2r}\int_t^T e^{i(t-s)(P_{k}-V)}\chi(x)\psi(D_x/k)f(x,s)ds,$$
so that
$$(D_t+P_{k}-V)v_j=\pm i k^{2r} \chi\psi f,$$
and it suffices to estimate
$$\|\psi \chi v_j\|_{L^2L^2}\leq C\| f\|_{L^2L^2}.$$
Now taking the Fourier transform in time and using Plancheral's theorem, we have that it suffices to estimate
$$\|\psi \chi \hat{v_j}\|_{L^2L^2}\leq C \|\hat{f}\|_{L^2L^2}$$
but this is the same as estimating
$$\|\psi\chi k^{2r}(\tau \pm i0 +P_{k}-V)^{-1}\chi\psi\|_{L^2_x\ra L^2_x}\leq C.$$
This means that for the operator $P_k$ defined above we can reduce the estimate to showing that
$$\|\psi \chi k^{2r}(\tau \pm i0 + P_k-V)^{-1} \chi \psi \|_{L^2_x\ra L^2_x}\leq C.$$ 
Let $-z = \tau k^{-2}$ and $h=k^{-1}$ to get
$$\|\psi \chi (-z \pm i0 + (hD_x)^2+V_1-h^2V_2\p_\omega+h^2V)^{-1} \chi \psi \|_{L^2_x\ra L^2_x}\leq C.$$

In particular we want to show that
\begin{equation}
  \|(-z+ (hD_x)^2+V_1-h^2V_2\p_\omega+h^2V)\chi(x)\psi(hD_x)\varphi
  u\|^2_{L^2} \geq h^{\frac{4m+2}{2m+3}}\|\chi(x)\psi(hD_x)
  (\phi u)\|^2_{L^2}. \label{E:ml-1}
  \end{equation}
So with the following lemma we can get the desired result.
\begin{lemma}
  \label{L:ChMe}
For $\varepsilon>0$ sufficiently small, let $\varphi\in S(T^*\IR)$ have compact support in $\{|(x-1,\xi)\leq \varepsilon\}$. Then there exists $C_\varepsilon>0$ such that
\begin{equation}
  \label{E:ChMe}
  \|(P-z)\varphi^w u\|^2 \geq C_\varepsilon
  h^{(4m+2)/(2m+3)}\|\varphi^wu\|^2,\,\,z\in [C-\varepsilon,
    C+\varepsilon]
  \end{equation}
where $P=(hD_x)^2+V_1-h^2V_2\partial_\omega^2-h^2V$ and $\varphi^w$ denotes quantization in only the $x$ and $\p_x$ directions and $\|\cdot\|$ denotes the$ L^2$ norm in $x$ and $\omega$ coordinates.
\end{lemma}


Now when looking at the norm we absorb the $h^2V$ term to the right
hand side  of \eqref{E:ChMe} since $(4m +2)/(2m +3) <2$. We just
need to deal with the $-h^2V_2\partial_\omega^2$ term because
$$\|((hD_x)^2+V_1)\varphi^w u\|_{L^2}^2\geq
h^{(4m+2)/(2m+3)}\|\varphi^w u\|_{L^2}^2$$ by \cite{ChMe-lsm}.

Now we very briefly summarize  the commutator process as in
\cite{ChMe-lsm}. We define $\Lambda, \Lambda_2$ as follows: freeze $\epsilon_0>0$ and let 
\[
\Lambda(r) = \int_0^r \lll t \rrr^{-1-\epsilon_0} dt, \,\,\,
\Lambda_2(r) = \int_{-\infty}^r \lll t \rrr^{-1-\epsilon_0} dt.
\]
For the remainder of the paper, we denote by $\chi(s)$ a smooth, even
bump function with $\chi(s) \equiv 1$ for $|s| \leq \delta_1$ and
support in $\{ |s| \leq 2 \delta_1 \}$.  Here $\delta_1 \gg \epsilon$
where $\epsilon>0$ is as in Lemma \ref{L:ChMe}.
with compact support near $0$ so that $\chi(x-1) \chi(\xi)$
microlocalizes (in the semi-classical sense)
near $x = 1$ and $\xi = 0$.  
Just as in \cite{ChMe-lsm}, let $\th \gg h$ be a second small
parameter and let 
$$a(x,\xi; h, \tilde{h})=\Lambda(\Xi)\Lambda_2(X-1)\chi(x-1)\chi(\xi)$$
where
$$X-1=\frac{x-1}{(h/\tilde{h})^\alpha},\,\, \Xi=
\frac{\xi}{(h/\tilde{h})^\beta}.$$
Here
\[
\alpha = \frac{2}{2m_1 + 3}, \,\, \beta = 1-\alpha.
\]

We now employ a similar commutator method to get a favorable sign
on the $V_2$ term.  To somewhat ease notation, let $v = \phi^w u$.  We
have that \begin{equation}
  \label{E:V2-est-1}
  C\|v\| \|Pv\|\geq \<[P,a^w],v,v\> =
  \<i[(hD_x)^2+V_1,a^w]v,v\>+\<i[-h^2V_2,a^w]\p_\omega^2 v,  v\>
  \end{equation}
for some constant $C$.  The first term is exactly the same as in
\cite{ChMe-lsm}.  
We have assumed that $V_2$ is decreasing and linear near $x = 1$, so
if $\epsilon>0$ is sufficiently small, $V_2''$ is supported away from
the support of $\chi(x-1)$.  
In particular, we have $i[-V_2,a^w]= -h(H_{V_2}a)^w$,
where
\[
- H_{V_2} a = V_2'(x) \p_\xi a \leq 0
\]
on the wavefront set of $v$.  Hence
\[
\lll (-H_{V_2}a)^w \p_\omega^2 v , v \rrr = \lll (H_{V_2} a)^w
\p_\omega v , \p_\omega v \rrr \geq 0.
\]
Plugging in to \eqref{E:V2-est-1}, 
this implies that
$$C\|v\| \|Pv\|\geq  \<i[P,a^w]v,v\> \geq \<i[hD_x^2+V_1,a^w]v,v\>.
$$
Hence, by the results in \cite{ChMe-lsm} we have that
$$\|v\|\, \|Pv\|\geq \frac{1}{C}h^{(4m+2)/(2m+3)}\|v\|^2+\frac{C'}{C}h^{3-\beta}\|\p_\omega v\|^2$$

This completes the proof when separating variables in the $\theta $
direction.  Separating variables in just the $\omega $ direction is similar.



\bibliographystyle{alpha}
\bibliography{MW-bib}

\end{document}